\documentclass[draft]{amsart}
\usepackage{amssymb, enumerate, xspace, graphicx,stmaryrd}
\usepackage[all]{xy}
\usepackage[english]{babel}
\usepackage[T1]{fontenc}
\usepackage[latin1]{inputenc}
\SelectTips{cm}{}

\numberwithin{equation}{section}

1

\numberwithin{subsection}{section}

\allowdisplaybreaks[1]

\newtheorem*{namedtheorem}{\theoremname}
\newcommand{\theoremname}{testing}

\newtheorem{theorem}{Theorem}[section]

\newtheorem{proposition}[theorem]{Proposition}
\newtheorem{proposition-definition}[theorem]
{Proposition-Definition}
\newtheorem{corollary}[theorem]{Corollary}
\newtheorem{lemma}[theorem]{Lemma}

\theoremstyle{definition}
\newtheorem{definition}[theorem]{Definition}

\newtheorem{remark}[theorem]{Remark}

\theoremstyle{remark}

\newcommand\cA{\mathcal{A}} 
\newcommand\cC{\mathcal{C}} \newcommand\cD{\mathcal{D}}
 \newcommand\cF{\mathcal{F}}
 
\newcommand\cI{\mathcal{I}} 
 \newcommand\cL{\mathcal{L}}
\newcommand\cM{\mathcal{M}} 
\newcommand\cO{\mathcal{O}} 
 
\newcommand\cS{\mathcal{S}} \newcommand\cT{\mathcal{T}}
\newcommand\cU{\mathcal{U}}

 \newcommand\DD{\mathbb{D}}

 \newcommand\PP{\mathbb{P}}

 \newcommand\ZZ{\mathbb{Z}}

\newcommand\OO{\mathcal{O}_C}

\newcommand\arr{\ifinner\to\else\longrightarrow\fi}

\newcommand\arrto{\ifinner\mapsto\else\longmapsto\fi}

\def\displaytimes_#1{\mathrel{\mathop{\times}\limits_{#1}}}

\def\displayotimes_#1{\mathrel{\mathop{\bigotimes}\limits_{#1}}}

\newcommand\pic{\operatorname{Pic}}

\newdir{ >}{{}*!/-5pt/@{>}}

\newcommand\doublelong[2]{\mathbin{\xymatrix{{}\ar@<3pt>[r]^{#1}
\ar@<-3pt>[r]_{#2}&}}}

\newlength{\ignora}

\newcommand{\lra}{\longrightarrow}

\newcommand{\ra}{\rightarrow}

\newcommand{\uu}{\mathcal{U}_C(r,0)} 

\newcommand{\su}{\mathcal{SU}_C(r)}

\newcommand{\urg}{\mathcal{U}_C(r,rg)}

\newfont{\sheaf}{eusm10 scaled\magstep1}

\begin{document}

\title[Modular and birational geometry of $\mathcal{SU}_C(r)$]{Modular subvarieties and birational geometry of $\mathcal{SU}_C(r)$.}

\author[Bolognesi]{Michele Bolognesi}
\address{Institut f\"{u}r Mathematik\\Humboldt Universit\"{a}t zu Berlin\\Rudower Chaussee 25\\10099 Berlin\\Germany\\ Università Roma Tre\\ Dipartimento di Matematica\\ L.go S. Leonardo Murialdo 1\\00146 Roma\\ Italy
}
\email[Bolognesi]{michele.bolognesi@sns.it}

\author[Brivio]{Sonia Brivio}
\address{Dipartimento di Matematica "F.Casorati"\\Università di Pavia\\Via Ferrata 1\\
27100\\ Italy}
\email[Brivio]{sonia.brivio@unipv.it}


\begin{abstract}
Let $C$ be an algebraic smooth complex genus $g>1$ curve. The object of this paper is the study of the birational structure of the coarse moduli space $\mathcal{U}_C(r,0)$ of semi-stable rank r vector bundles on $C$ with degree 0 determinant and of its moduli subspace $\mathcal{SU}_C(r)$ given by the vector bundles with trivial determinant. Notably we prove that $\uu$ (resp. $\su$) is birational to a fibration over the symmetric product $C^{(rg)}$ (resp. over $\PP^{(r-1)g}$) whose fibres are GIT quotients $(\PP^{r-1})^{rg}//PGL(r)$.  In the cases of low rank and genus our construction produces families of classical modular varieties contained in the Coble hypersurfaces.
\end{abstract}

\maketitle

\section{Introduction}

Let $C$ be a genus $g>1$ smooth complex algebraic curve, if $g\neq 2$ we will also assume that $C$ is non-hyperelliptic. Let $\uu$ be the moduli space of rank r semi-stable vector bundles on $C$ with degree zero determinant and let us denote as usual $\su$ the moduli subspace given by vector bundles with trivial determinant. These moduli spaces appeared first in the second half of the last century thanks to the foundational works of 
Narashiman-Ramanan  \cite{ramanara} and Mumford-Newstead \cite{MumNews} and very often their study has gone along the study of the famous theta-map

\begin{eqnarray*}
\theta:\su & \dashrightarrow & |r\Theta|;\\
E & \mapsto & \Theta_E:=\{L\in \pic^{g-1}|h^0(C,E\otimes L)\neq 0\}.
\end{eqnarray*}

While we know quite a good deal about $\theta$ for low genuses and ranks, as the rank or the genus grow our knowledge decreases dramatically, see Sect. 2 for a complete picture of known results.
The question of rationality is even more haunting. When rank and degree are coprime the situation is quite settled \cite{daking} but when the degree is zero (or degree and rank are not coprime) the open problems are still quite numerous. It is known that all the spaces $\su$ are unirational but the rationality is clear only for $r=2,\ g=2$, when in fact the moduli space is isomorphic to $\PP^3$, \cite{ramanara}. Some first good ideas about the birational structure for $g=2$ were developed in \cite{ang}. Then the $r=2$ case was analyzed in any genus by the first named author and A.Alzati in \cite{albol} with the help of polynomial maps classifying extensions in the spirit of \cite{ab:rk2}. In this paper we give a description for the higher rank cases.

\begin{theorem}\label{global}
Let $C$ be a smooth complex curve of genus $g>1$, non-hyperelliptic if $g>2$, then $\uu$ (resp. $\su$) is birational to a fibration over $C^{(rg)}$ (resp. $\PP^{(r-1)g}$) whose fibers are GIT quotients $(\PP^{r-1})^{rg}//PGL(r)$.
\end{theorem}



In the case of $\su$, Proposition \ref{global} allows us to give a more precise explicit description of the fibration of $\su$ in the case $r=3,g=2$. In fact $\mathcal{SU}_C(3)$ is a double covering of $\PP^8$ branched along a hypersurface of degree six $\cC_6$ called the \textit{Coble-Dolgachev sextic} \cite{ivolocale}. Our result is the following.

\begin{theorem}
The Coble-Dolgachev sextic $\cC_6$ is birational to a fibration over $\PP^4$ whose fibers are Igusa quartics.
\end{theorem}

We recall that an Igusa quartic is a modular quartic hypersurface in $\PP^4$ that is related to some classical GIT quotients (see e.g. \cite{do:pstf}) and moduli spaces. Its dual variety is a cubic 3-fold called \textit{Segre cubic}, that is isomorphic to the GIT quotient $(\PP^1)^6//PGL(2)$. 

\smallskip

If $r=2$ and $g=3$, then $\cS\cU_C(2)$ is embedded by $\theta$ in $\PP^7$ as a remarkable quartic hypersurface $\cC_4$ called the \textit{Coble quartic}, \cite{ramanara}. Our methods also allow us to give a quick proof of the following fact, already shown in \cite{albol} by means of polynomial maps.

\begin{proposition}
The Coble quartic $\cC_4$ is birational to a fibration over $\PP^3$ whose fibers are Segre cubics.
\end{proposition}

We underline that the cases of $\cC_4$ and $\cC_6$ are particularly interesting because one can interpret the beautiful projective geometry of the Igusa quartic and the Segre cubic in terms of vector bundles on $C$ (see Sect. 6). We hope that these results could help to shed some new light on the question of rationality of $\su$ and on the properties of the theta map.

\smallskip

\textit{Acknowledgments:} The inspiration of this work came from a conversation that the first named author had with Norbert Hoffmann. This work also enjoyed from remarks and discussions with Christian Pauly, Alessandro Verra and Angela Ortega. Arnaud Beauville brought to our attention some existing literature. Quang Minh and Cristian Anghel kindly let us have copies of their papers. Finally Duco Van Straten and Edoardo Sernesi gave us good advices about the deformation theory of projective singular hypersurfaces.

\medskip

\textbf{Description of contents.}

In section 2 we collect a few results on the theta map. In section 3 we outline the relation between theta maps and theta-linear systems by introducing the theta divisor of a vector bundle with integral slope. 
 In section 4 we  describe some properties of the evaluation and the determinant map for a vector bundle of any rank. In Section 5 we prove Theorem \ref{global} and discuss briefly the relation between slope stability and GIT stability. Finally, in Section 6 we apply the results to the cases $g=2,\ r=3$ and $g=3,\ r=2$ and give an explicit description of the fibration. 

\section{The Theta map.}\label{thetas}

Let $C$ be a smooth complex  algebraic curve of genus $g \geq 2$, we assume  that it is  non-hyperelliptic if $g>2$. 
Let $\pic^d(C)$ be the Picard variety  parametrizing  line bundles of  degree $d$  on $C$,  $\pic^0(C)$ will be often denoted as $J(C)$. Let  $\Theta \subset \pic^{g-1}(C)$ be the canonical theta divisor   

$$\Theta:=\{L\in \pic^{g-1}(C)| \ h^0(C,L)\neq 0 \}.$$

For $r \geq 2$, let $\su$ denote the coarse moduli space of semi-stable vector bundles of rank $r$ and trivial determinant on $C$. It is a normal, projective variety  of dimension $(r^2-1)(g-1)$. 
It is well known that $\su$ is locally factorial and that $\pic (\su) = \ZZ$ \cite{dn:pfv}, generated by a line bundle
$\cL$ called the \textit{determinant bundle}. On the other hand, for $E\in \su$ we define

$$\Theta_E:=\{L\in \pic^{g-1}(C) | \ \ \ h^0(C,E\otimes L)\neq 0 \ \ \}.$$

This is either a divisor in the linear system $|r\Theta|$ or the whole $\pic^{g-1}(C)$. For $E$ a general  bundle  $\Theta_E$ is a divisor, {\it the theta divisor} of $E$.  This means that one can define the rational {\it theta map} of $\su$:

\begin{eqnarray}
\theta: \su  \dashrightarrow  |r\Theta|
\end{eqnarray}

sending $E $ to its theta divisor  $\Theta_E$.
The relation between the theta map and the determinant bundle  is given by the following fundamental result:

\begin{theorem}\cite{bnr}
There is a canonical isomorphism $|r\Theta| \stackrel{\sim}{\lra} |\cL|^*$ which  identifies $\theta$ with the rational map ${\varphi}_{\cL} \colon 
\su \dashrightarrow |\cL|^*$ associated to the determinant line bundle. 
\end{theorem}


The  cases when $\theta$ is a morphism or finite are of course very appealing. 
Notably,  $\theta$ is an embedding for $r=2$ \cite {ramanara}, \cite{BVtheta},\cite{vgiz} and it is a morphism when $r=3$ for $g \leq 3$ and  for a general curve of genus $g > 3$, \cite{bove23}, \cite{raysect}. Finally, $\theta$ is generically finite for $g=2$  \cite{bove23}, \cite{briverraBN} and we know its  degree  for $r \leq 4$,  \cite{ivolocale}, \cite{deg16}. 
There are also good descriptions of the image of $\theta$ for $r=2 \quad g=2,3$ \cite{rana:cra} \cite{paulydual},  $r=3,\ g=2$ \cite{orte:cob} \cite{Quang},  $r=2,\ g=4$ \cite{oxpau:heis}.  
Moreover, it has recently been shown in \cite{briverraPL} that if $C$ is general and $ g >> r$ then $\theta$ is generically injective.

\section{Vector bundles and theta linear systems}\label{thetasystems}

The notion of theta divisor can be extended to  vector bundles  with integral slope.
Let $\urg$ be the  moduli space of semi-stable vector bundles on $C$ with rank $r$ and degree $rg$.
The tensor product defines a natural map:

\begin{eqnarray*}
t \colon  \su \times \pic^{g}(C) & \ra & \urg;\\
(E, \cO_C(D)) & \mapsto  &  E \otimes {\cO}_C(D), 
\end{eqnarray*}

which is étale, Galois, with Galois group $J(C)[r]$, the group of $r$-torsion points of the Jacobian of $C$. 

   
Moreover, if one restricts $t$ to $\su \times \cO_C(D)$ this yields an isomorphism $t_{D}  \colon \su \to  \mathcal{SU}_C(r,\cO_C(rD)),$ where the latter is the moduli space of rank $r$ semi-stable vector bundles with determinant  $\cO_C(rD)$.

\begin{definition}\label{supertheta}
Let $F \in \urg$, then we define the {\it theta divisor } of $F$ as follows:

\begin{equation*}
\Theta_F:=\{L \in \pic^1(C) | \ h^0(C,F \otimes L^{-1} )\neq 0\}.
\end{equation*} 

\end{definition}

Let $E \in \su$ and $\cO_C(D) \in \pic^g(C)$. If $F = E \otimes \cO_C(D)$, then we have that $\Theta_F = \cO_C(D) - {\Theta}_E$, thus $\Theta_F$ is a divisor if and only  ${\Theta}_E$ is a divisor.  We define 

\begin{equation} 
{\Theta}_D \colon = \{ L \in \pic^1(C) \vert \ \ h^0(\cO_C(D) \otimes L^{-1}) \geq 1 \ \}.
\end{equation}

Then, for any $r \geq 1$,  we have a natural isomorphism $\sigma_D:|r\Theta| \ra |r\Theta_D|$ given by the translation
$M  \mapsto \cO_C(D) - M$. Moreover, if $\cO_C(rD_1) \simeq \cO_C(rD_2)$, then $|r\Theta_{D_1}| = |r\Theta_{D_2}|$ so we conclude that if $F  \in \urg$ admits theta divisor, then ${\Theta}_F  \in \vert r {\Theta}_D \vert$, for any line bundle $\cO_C(D) \in  \pic^g(C)$ which is a $r$-root of $\det F$. 
In this way we obtain a family of theta linear systems over the Picard variety $\pic^{rg}(C)$, as the following shows.




\begin{lemma} 
There exists a projective bundle ${\mathcal T}$  over $\pic^{rg}(C)$:
$$ p \colon {\mathcal T} \to \pic^{rg}(C), $$
whose fiber over $\cO_C(M)\in \pic^{rg}(C)$  is  the linear system $ \vert r {\Theta}_D \vert$, 
where $\cO_C(D)\in \pic^{g}(C) $ is any $r$-root of $\cO_C(M)$.

\end{lemma}

\begin{proof}
Remark first that the linear system  $|r\Theta_D|$ is well-defined since it does not depend on which $r$-root $\OO(D)$ of $\OO(M)$ we choose. Let us now consider the tensor product map:

\begin{eqnarray*}
\delta \colon \pic^g(C) \times \pic^1(C) & \to & \pic^{g-1}(C);\\
(\cO_C(D), L) & \mapsto & \cO_C(D) \otimes L^{-1}.
\end{eqnarray*}

For any $\cO_C(D) \in \pic^g(C)$,  we have: ${\delta}^* {r\Theta}_{|\cO_C(D) \times \pic^1(C)} \simeq \cO_{\pic^1(C)} ({r\Theta}_D).$ Let $p_1 \colon  \pic^g(C) \times \pic^1(C) \to \pic^g(C)$ be the projection onto the first factor. Consider the sheaf $ \cF:={p_1}_{*} \cO_{\pic^g(C) \times \pic^1(C)}( {\delta}^*(r\Theta)).$ It is locally free and its fiber at $\cO_C(D) \in \pic^g(C)$ is canonically identified with  the vector space $H^0( \pic^1(C), \cO_{\pic^1(C)} (r{\Theta}_D))$. We call $\widetilde{{\mathcal T}}$ the projective bundle $\PP (\cF)$ on $\pic^g(C)$. 

Moreover the vector bundle $\cF$ is $J[r](C)$-equivariant, hence by easy descent theory (see \cite{stacchi} Thm. 4.46) it passes to the quotient by $J[r](C)$, i.e. the image of the cover $\rho  \colon \pic^g(C) \to \pic^{rg}(C)$ given by taking the $r^{th}$ power of each $L\in \pic^g(C)$. The projectivized of the obtained bundle is the projective bundle $\cT$ we are looking for.
We denote by  $p:\cT \ra \pic^{rg}(C)$  the natural projection on the base of the projective bundle.


\end{proof}

The previous arguments allow us to define the rational {\it theta map} of $\urg$.

\begin{eqnarray}
 \theta_{rg} \colon \urg &  \dashrightarrow  &  {\mathcal T};\\
 F & \mapsto & \Theta_F.
\end{eqnarray}


Let us denote by $\theta_D$ the restriction of $\theta_{rg}$ to the  subspace $\mathcal{SU}_C(r,\cO_C(rD))$. Then we have the following commutative diagram:

$$\xymatrix{ \su \ar[r]^(.4){t_{D}} \ar@{-->}[d]_{\theta} & \mathcal{SU}_C(r,\cO_C(rD))  \ar@{-->}[d]^{{\theta}_D} \\
\vert r \Theta \vert \ar[r]^{{\sigma}_{D}} & \vert r {\Theta}_{D} \vert \\}$$

since $t_{D}$ and ${\sigma}_{D}$ are isomorphism, we can identify the two theta maps. 
Finally remark that the composed  map $p \circ \theta_{rg}$ is precisely the natural map $det \colon  \urg \to \pic^{rg}(C),$ which associates to each vector bundle $F$ its determinant line bundle $det (F)$.

\section{The fundamental divisor of a vector bundle}

In the sequel, where we don't state differently, we will be concentrating on semi-stable vector bundles on $C$ with rank $r$ and degree $rg$.
 Let  $F \in \urg$, note that $\chi(F) = r$ hence  $h^0(F)  \geq  r$.  We  can associate to $F$  two natural maps. The first one is  the {\it evaluation  map} of $F$:

\begin{eqnarray*}
ev_F \colon H^0(F)  \otimes \cO_C & \lra & F;\\
(s,x) & \mapsto & s(x).
\end{eqnarray*}

\noindent The second one is  the {\it determinant map} of $F$:

\begin{eqnarray*}
d_F \colon {\wedge}^rH^0(F) & \lra & H^0(det F);\\
s_1 \wedge s_2 ... \wedge s_r & \mapsto & ( x \mapsto s_1(x)  \wedge s_2(x)  \wedge... \wedge s_r(x) ).
\end{eqnarray*}

The image of the determinant map defines  a linear subsystem  $|Im(d_F)| \subset |det(F)|$.
These two maps are somehow dual to each other and some properties of the evaluation map can be translated in the language of the determinant map and vice versa, as the following Lemma shows.

\begin{lemma}\label{dege}

Let $F \in \urg$, then we have:

\begin{enumerate}
\item{} $rk{(ev_F)}_x \leq r-1$ for a point $x \in C$ if and only if
$x$ is a base point of $|Im (d_F)|$;

 

\item $rk{(ev_F)} \leq r-1$ if and only if $d_F$ is the zero map;


\item{} if $h^0(F) = r$, then  either $d_F$ is the zero map or  $\vert Im (d_F) \vert =   \{  \sigma \} $, for a non zero section  $\sigma \in |detF|$. In this case  $ev_F$ is  generically surjective  and its degeneracy locus is the zero locus of $\sigma$.
\end{enumerate}

\end{lemma} 

\begin{remark}\label{g2}
 Let  $g= 2$ and   $F \in {\mathcal U}_C(r,2r)$. The condition $h^0(F) > r$ is equivalent to 
$Hom(F,{\omega}_C) \not= 0$. Since $\mu(F) = \mu ({\omega}_C)$ this implies that $F$ is not stable and $F = {\omega}_C \oplus G$, for some $G \in {\mathcal U}_C(r-1,2(r-1))$.
\end{remark}

\noindent Let us now define the following subset of $\urg$:

$$ \mathcal U = \{ F \in \urg \ \vert \ h^0(F) = r \ \ \  d_F \not\equiv 0 \ \}.$$

\noindent First we show that $\cU \subset \urg$ is not empty.

\begin{lemma}
$\cU$ is a non-empty open subset of $\urg$.
\end{lemma}

\begin{proof}
Since the conditions that define $\cU$ are open, it is enough to produce a semi-stable vector bundle $F \in \cU$.
For example, let $F = L_1 \oplus... \oplus L_r$ be the sum of $r$ line bundles of degree $g$ with $h^0(L_i) = 1$ for any $i = 1,.,r$. Then we have $H^0(F) = {\bigoplus}_{i = 1}^{r} H^0(L_i)$ and $h^0(F) = r$. 
Moreover $d_F$ is just the natural multiplication map of global sections

\begin{eqnarray*}
\nu \colon H^0(L_1) \otimes H^0(L_2) \otimes .. H^0(L_r) & \to & H^0(L_1 \otimes.. \otimes L_r);\\
(s_1,\dots, s_r) & \mapsto & s_1\cdot \dots \cdot s_r;
\end{eqnarray*}

\noindent that is non-zero.
\end{proof} 

Now suppose that we have a vector bundle $F \in \cU$. 
Let  $\{ s_1, s_2,..,s_r \}$ be a base of $H^0(F)$ and  let  $\sigma = d_F(s_1 \wedge... \wedge s_r)$. Then Lemma \ref{dege} implies that the divisor ${\mathbb D}_{F}:= \mathrm{Zeros(\sigma)} \in |\det F|$ is well defined.  





\begin{definition}
We call ${\mathbb D}_{F}$ the {\it fundamental divisor} of  $F \in {\mathcal U}$.
\end{definition}

Remark that if $F\in \urg$ than it has integral slope and we can consider the theta divisor ${\Theta}_{F} = \{ L \in \pic^{1}(C) \ \vert h^0(C, F \otimes L^{-1}) \not= 0 \}$ of Def. \ref{supertheta}. 
Now the natural Abel-Jacobi map $a \colon C \ra  Pic^1(C)$ embeds $C$ in $\pic^1(C)$, and with a little argument  one can show the following.



\begin{proposition}\label{pulbec}

Let $F \in \urg$. Then $F \in  \cU$  if and only if  $F$ admits theta divisor ${\Theta}_{F}$ and $a(C) \not\subset {\Theta}_{F}$. In this case we have:

$$\mathbb D_{F} = a^* ( {\Theta}_{F}).$$

\end{proposition}
 






\section{The  fundamental map and its fibers}\label{secdet}

\newcommand{\xrg}{C^{(rg)}}

Let $C^{(rg)}$ be the $rg$- symmetric product of $C$. We recall 
\cite{acgh:gac} that for every $r\geq 2$, $\xrg$ has a natural structure of  projective bundle over  $\pic^{rg}(C)$, given by the natural Abel map $a_{rg} \colon \xrg  \ra   \pic^{rg}(C)$. The fiber over $\cO_C(M) \in \pic^{rg}(C)$ is the complete linear system $\vert \cO_C(M) \vert$.



\begin{definition}
We call 

\begin{eqnarray*}
\Phi \colon \cU &  \lra & C^{(rg)},\\ 
F & \mapsto & \mathbb{D}_{F}
\end{eqnarray*}

the {\it fundamental map} of $\urg$.
\end{definition}

\noindent The aim of this section is the description of the fibers of $\Phi$. We start by showing some basic properties of the map $\Phi$ itself. First of all, note that since $\mathbb{D}_{F} \in \vert det F \vert$, then we have a commutative diagram:

$$\xymatrix{ \urg \ar[r]^{\Phi} \ar[d]_{det} &  \xrg  \ar[d]^{a_{rg}} \\
\pic^{rg}(C) \ar[r]^{id} & \pic^{rg}(C). }$$

Let $D\in \pic^g(C)$, then the restriction of $\Phi$ to the moduli space $\mathcal{SU}_C(r,\cO_C(rD))$ induces a map ${\Phi}_D \colon \mathcal{SU}_C(r,\cO_C(rD))  \dashrightarrow  \vert \cO_C(rD) \vert$. 

\begin{theorem}\label{domino}
For any $r \geq 2$ and $g \geq 2$, $\Phi$ is dominant. For any $\cO_C(D) \in \pic^g(C)$, $\Phi_D$ is defined and dominant. 
\end{theorem}

\begin{proof} 
Remark that if $G \in C^{(rg)}$ can be written as the sum of $r$ non special effective divisor $G_i$ of degree $g$, $i=1,\dots, r$, then $G = {\Phi}(F)$,  where    $F = \bigoplus_{i=1}^r \cO_C(G_i)$. So the assertion follows  once we prove that for any $r \geq 1$ a general  divisor of  $C^{(rg)}$ satisfies the above property.
We prove, by induction on $r$, that   the set $S_{rg} \subset C^{(rg)}$ of divisors which does not  satisfy this property 
is a proper closed subvariety. 
\hfill\par\noindent
Let $r= 1$, then $S_g \subset C^{(g)}$ is the closed set of special divisor of degree $g$ and $\dim  S_g = g-2$ \cite{acgh:gac}. For any $r \geq 2$, let us consider the sum map: 
$$m \colon C^{(g)} \times C^{((r-1)g)} \to C^{(rg)}, \quad (D_1, D_2) \to D_1 + D_2.$$
 This is a surjective finite map and we have:
$$S_{rg} = m(C^{(g)} \times S_{(r-1)g}) \cup m(S_g \times C^{((r-1)g)}).$$ 
By induction hypothesis $ S_{(r-1)g} $ is a proper closed subvariety of $C^{((r-1)g)}$ and $S_g$ is of codimension 2, hence $S_{rg}$ is a proper closed subvariety of $C^{(g)}$.
Finally, we prove that for $r \geq 2$ and for any $M \in  \pic^{rg}(C)$, the closed variety $S_{rg}$ does not contain all the elements of the linear system $\vert M \vert$. Let us consider
 $$m^*(\vert M \vert) = \{ (D_1, D_2) \in C^{(g)} \times C^{((r-1)g)} \vert \ \ D_1 + D_2 \in  \vert M \vert \}.$$
By Riemann-Roch the natural projection $p_2 \colon m^*(\vert M \vert) \to C^{((r-1)g)}$ is surjective and  $\dim m^*(\vert M \vert) =  \dim C^{((r-1)g)}$, thus $p_2$ is generically finite. Since the fiber of $p_2$ at $D_2 \in C^{((r-1)g)}$ is the linear system $\vert M(-D_2) \vert$, then it is connected. This in turn implies that for a general $D_2$, $p_2^{-1}(D_ 2) = \{ D_1 \}$,  hence   $D_1$ is not special.  This implies the claim. 
\end{proof}

Let now $F \in \mathcal{SU}_C(r,\cO_C(rD))$ be a general bundle admitting theta divisor $\Theta_F$. By Proposition \ref{pulbec}, the pull back defines a linear projection map 

\begin{equation}\label{proietta}
 a^* : |r{\Theta}_D| \to    \vert \cO_C(r D) \vert,
\end{equation} 

and the map $\Phi_D$ is the composition of $a^*$ with the theta map on $\mathcal{SU}_C(r,\cO_C(rD))$. The projection of Eq. \ref{proietta} can be globalized in a rational map $A^* \colon {\mathcal T} \to \xrg$ between projective bundles over $\pic^{rg}(C)$ that restricts to the corresponding projection on each fiber. This implies the following result.



\begin{proposition} 
The map ${\Phi}: \urg \dashrightarrow  C^{(rg)} $ is the composition of ${\theta}_{rg}$ with $A^*$. 
In particular, the map ${\Phi}_D: \mathcal{SU}_C(r,\cO_C(rD)) \dashrightarrow \vert \cO_C(rD) \vert$ 
is the composition of $\theta_D$ 
with  $a^*$.
\end{proposition}

As usual, let $({\PP}^{r-1})^{rg} // PGL(r)$ denote  the GIT quotient of $({\PP}^{r-1})^{rg}$ by the diagonal action of $PGL(r)$. We recall,  (\cite{do:pstf}, Thm. 1 pag. 23), that a  point $v= (v_1,..,v_{rg})  \in ({\PP}^{r-1})^{rg}$   is GIT semi-stable (resp. stable) if and only if for any subset $\{v_1,\dots,v_k\}$ of $v$ we have $ dim(Span(v_1,\dots,v_k)) \geq \frac{k}{g}$ (resp. >).  Then we have the following result:

\begin{theorem}\label{genfib}
The general fiber of $\Phi$ is birational to $({\PP}^{r-1})^{rg} // PGL(r)$.
\end{theorem}

Let us remark the fact that in general we have only a birationality result, whereas (see Sect. 6) if $g,r\leq 3$ and $g\neq r=3$ the general fiber is \textit{isomorphic} to the corresponding GIT quotient.


\begin{proof}
Let  $B \in C^{(rg)}$ be a general divisor in the image of $\Phi$.   For  simplicity, we assume  $B$ out of the big diagonal $\Delta$, that is 
$B = \sum_{i=1}^{rg}{x_i}$, with $x_i \not= x_j$, then:
$${\Phi}^{-1}(B) =  \{F \in {\mathcal U}  \colon {\mathbb D}_F = B \}.$$ 
Since $\Phi$ is dominant,we have:

$$dim {\Phi}^{-1}(B)= dim (\urg)  - dim (C^{(rg)})= (r^2-r)g - (r^2 -1),$$

which is actually the dimension of the variety  $({\PP}^{r-1})^{rg} // PGL(r)$.
Let  $F \in {\Phi}^{-1}(B)$, since $B$ is the degeneracy locus of the evaluation map, by dualizing we find the following exact sequence.

\begin{equation}\label{sequa2}
0 \to F^* \to H^0(F)^* \otimes \cO_C \to \cO_B \to 0.
\end{equation}

This means that, up to the choice of a basis  of $H^0(F)$, $F^*$ is the kernel of a surjective morphism in $\mathrm{Hom}(\cO_C^{ r},\cO_B).$
Let us now consider a vector space $V$ of dimension $r$ and the   projective space  $\PP(V^*)$. By mimicking sequence \ref{sequa2}, we can construct a flat family of vector bundles on $C$ over $(\PP(V^*))^{rg}//PGL(r)$ in  the following way.
 We take $v=(v_1,\dots,v_{rg}) \in \PP(V^*)^{rg}$ and  let $\phi_v$ be the surjective morphism of sheaves $V\otimes \cO_C \ra \cO_B$ that is the zero map out of the support of $B$ and that is obtained by taking one lift of $v_i$ to $V^*$ and applying it on the fiber of $V\otimes \cO_C$ over $x_i \in B$. The morphism $\phi_v$ depends on the choice of the lift but the kernel of the sequence

\begin{equation}\label{univ}
0\lra Ker(\phi_v) \lra V\otimes \cO_C \stackrel{\phi_v}{\lra} \cO_B \to 0,
\end{equation}

is well defined over $\PP(V^*)^{rg}$. This implies that $E_v:=Ker(\phi_v)^*$, for  $v \in {\PP(V^*)}^{rg}$, is  a family of   rank $r$ vector bundles on $C$ with determinant $\cO_C(B)$. Moreover, it is invariant under the diagonal action of $PGL(r)$ on ${\PP(V^*)}^{rg}$. We will abuse slighlty notation by calling $E_v$ both the family of bundles over ${\PP(V^*)}^{rg}$ and, later, over its GIT quotient.


Remark that by definition,  any semi-stable vector bundle in the fiber $\Phi^{-1}(B)$ can be written as $E_v$,  for some set $v \in {\PP(V^*)}^{rg}$. We denote $A_b\subset \PP(V^*)^{rg}$ the open subset given by $v\in{\PP(V^*)}^{rg}$ s.t.  $E_v\in\Phi^{-1}(B)$. Moreover, recall that $(\PP^{r-1})^{rg}//PGL(r)$ is the quotient of the open semi-stable subset of $\PP(V^*)^{rg}$. We denote $A_s\subset \PP(V^*)^{rg}$ this open subset. Then, by dimensional reasons, $A_b\cap A_s \neq \emptyset$. By passing to the quotient (recall that the construction of $E_v$ is $PGL(r)$-invariant) this implies that there exists at least one semi-stable $E_v$, for $v \in {(\PP^{r-1})}^{rg}//PGL(r)$. In turn this implies, by the openness of semistability, that  
the generic element of the family  $E_v$, $v \in {(\PP^{r-1})}^{rg}//PGL(r)$, is a semi-stable vector bundle with  fundamental divisor $B$. By the universal property of the coarse moduli space,  this induces a birational map 

$$ (\PP^{r-1})^{rg}//PGL(r) \dashrightarrow   \urg.$$

\noindent which is regular and one to one on the quotient of the open set $A_b\cap A_s$.

\end{proof}


\noindent On the other hand, let us associate to $E\in \cU$ a point set $v_E\in (\PP(V^*))^{rg}//PGL(r)$ s.t., keeping the notation of the previous proof, $E^*$ is the kernel of $\phi_{v_E}$.

\begin{lemma}\label{stable}
If $E \in {\Phi}^{-1}(B)$  is stable (resp. semi-stable) then $v_E \in (\PP^{r-1})^{rg}$ is GIT stable (resp. semi-stable).
\end{lemma}

\begin{proof}
Suppose that $E$ is semi-stable and  $\DD_E=B$.  For  any subset $v':=\{v_1,\dots,v_k\}$ of $v_E$, let us denote $V_k^*:=Span(v_1,\dots,v_k)\subset V^*$ and $x_i$ the point of $B$ that correspond to $v_i$.  Then we get a commutative diagram:

$$\xymatrix{ 0  \ar[r] &  E^* \ar[d]   \ar[r]  & V \otimes \cO_C   \ar[d]    \ar[r]^{\phi_{v_{E}}} & \cO_B  \ar[d] \ar[r] & 0  \\
0  \ar[r]  &  G^*  \ar[r] \ar[d] &   V_k \otimes \cO_C   \ar[d] \ar[r]^{\phi_{v'}}  &   \oplus_{i=1}^k\cO_{x_i} \ar[d]  \ar[r] & 0 \\
& 0 & 0& 0 & \\ }$$ 

for some vector bundle $G^*$ with  $rk(G^*)=dim (V_k) = s$. Now,  in order to show that $v_{E}$ is GIT semi-stable, it is enough to show that $dim(V_k) \geq \frac{k}{g}$. Since we  have that $deg(G^*)= -k$ and $E^*$ is semi-stable, then $\mu(G^*) = \frac{-k}{s} \geq  -g$,  from which we obtain the desired inequality. The stable case is described in the same way but with strict inequalities.
\end{proof}

A part of Lemma \ref{stable} and Theorem \ref{genfib} is proved for $g=2$ in \cite{ang}. 
Finally, we describe the quotient of $A_s\cap A_b$ in $(\PP^{r-1})^{rg}//PGL(r)$. This is the sublocus of $(\PP^{r-1})^{rg}//PGL(r)$ where $E_v\in \Phi^{-1}(B)$. 
 
\begin{lemma} 
Let $v \in (\PP^{r-1})^{rg}//PGL(r)$,  $E_v \in {\phi}^{-1}(B)$   if and only if $h^0(E_v) = r$.
\end{lemma}

\begin{proof}
Let $v \in (\PP^{r-1})^{rg}//PGL(r)$. The only if part is clear. If $h^0(E_v) =r$, then by taking the dual of sequence \ref{univ}, we have that $H^0(E_v) \simeq V^*$,  the evaluation map of $E_v$ is generically surjective and its degeneracy locus is actually $B$. 
So it is enough to prove that $E_v$ is semi-stable. Suppose that there exists a  proper rank $s$ sub-bundle $F \subset E_v$ with $\mu(F)  > \mu(E_v)= g$. This by Riemann-Roch implies $h^0(F) \geq  s+1$. Since 
$H^0(F) \subset H^0(E_v)$  this contradicts  the  generically surjectivness of  the evaluation map of $E_v$.
\end{proof}

As a corollary we have the following result:
\begin{theorem}
The moduli space $\cU_C(r,0)$ (resp. ${\mathcal SU}_C(r)$) is birational to a fibration over $C^{(rg)}$ (resp. $\PP^{(r-1)g}$) whose fibers are $(\PP^{r-1})^{rg}//PGL(r)$.
\end{theorem}



\section{ Application to the Coble hypersurfaces}

\smallskip

When the genus of $C$ is 2 or 3 and the rank is small enough, the moduli spaces $\su$ have very nice explicit descriptions related to certain hypersurfaces, called the 
{\it Coble hypersurfaces}. 

The main theorems of this section show how these hypersurfaces are in fact fibrations over certain projective spaces whose fibers are isomorphic to classical modular varieties related to the moduli of 6 points on a line and on a plane: namely the {\it Segre Cubic} and the {\it Igusa quartic}.

\smallskip

The basic example of these fibrations, namely the $g=2,$ $r=2$ case, is quite instructive. See \cite{bol:kumwed} for details.




\subsection{The Coble Sextic}\label{sextic}

In this subsection we assume that  $C$ is a curve of genus $2$ and we consider the moduli space $\cS\cU_C(3)$  of semi-stable vector bundles on $C$ with  rank $3$ and trivial determinant. The theta map 

$$\theta:\cS\cU_C(3) \ra |3\Theta| \simeq \PP^8$$ 

is a finite  morphism   of  degree $2$ and the branch locus is a sextic hypersurface $\cC_6$, called Coble-Dolgachev sextic \cite{ivolocale}. The Jacobian variety $J(C)$ is embedded in $|3\Theta|^*$ as a degree 18 surface, and there exists a unique cubic hypersurface $\cC_3\subset |3\Theta|^*$ whose singular locus coincides with $J(C)$: {\it the Coble cubic}, \cite{coble}. 
It was then conjectured by Dolgachev, and subsequently proved in \cite{orte:cob} and indipendently in \cite{Quang}, that $\cC_6$ is the dual variety of $\cC_3$. 

\medskip

On the other hand the hypersurface known as Igusa quartic is the Satake compactification of the moduli space $\cA_2(2)$ of principaly polarized abelian surfaces with a level two structure \cite{ig:tc1}, embedded in $\PP^4$ by fourth powers of theta-constants. Anyway in our context appears rather because of its relation with the GIT geometry of sets of points in the projective plane. The GIT quotient of $(\PP^2)^6$ with respect to the diagonal action of $PGL(3)$ is a degree two covering of $\PP^4$ branched along a $\Sigma_6$-invariant quartic hypersurface $\cI_4 \subset \PP^4$, which is exactly the Igusa quartic. The involution that defines the covering is the Gale transform (also called association, for details see \cite{do:pstf}, \cite{eis:pop}), which is defined as follows. 

\begin{definition} (\cite{eis:pop}, Def. 1.1)

Let $r,s \in \ZZ$. Set $\gamma=r+s+2$, and let $\Gamma \subset \PP^r,$ $\Gamma' \subset \PP^s$ be ordered nondegenerate 
sets of $\gamma$ points represented by $\gamma \times (r+1)$ and $\gamma \times (s+1)$ matrices $G$ and $G'$, respectively. We say that $\Gamma'$ is the Gale transform of $\Gamma$ if there exists a nonsingular diagonal $\gamma \times \gamma$ matrix $D$ s.t. $G^T\cdot D \cdot G'=0$, where $G^T$ is the transposed matrix of $G$.

\end{definition}

The Gale transform acts trivially on the sets of 6 points in $\PP^2$ that lie on a smooth conic. The branch locus of the double covering is then, roughly speaking, the moduli space of 6 points on a conic and henceforth a birational model of the moduli space of 6 points on a line. One can say even more, in fact the GIT compactification of the moduli space of 6 points on a line is a cubic 3-fold in $\PP^4$, called the Segre cubic, and its projectively dual variety is the Igusa quartic (see \cite{koike}, \cite{hu:gsq} for details). From the projective geometry point of view the singular locus of $\cI_4$ is an abstract configuration of lines and points that make up a $15_3$ \textit{configuration}. This means the following: there are 15 distinguished lines and 15 distinguished points. Each line contains 3 of the points and by each point pass 3 lines (see \cite{dolgabstract} Sect. 9 for more). Moreover $\cI_4$ is the only hypersurface with such a singular locus in the pencil of $\Sigma_6$-invariant quartics in $\PP^4$ (\cite{huntnice}, Example 7).

\medskip

Let us twist as customary the vector bundles in $\cS\cU_C(3)$ by the degree 2 canonical bundle $K_C$ and get thus an isomorphism $\cS\cU_C(3)\cong \cS\cU(3,3K_C)$. Let $\Phi_K:\cS\cU(3,3K_C) \ra |3K_C|$ be the restriction of the fundamental map, then we have the following.



\begin{proposition}\label{pitre}
Let $L_e$ be the subsystem of theta divisors corresponding to decomposable bundles of type $E\oplus K_C$, with $E\in \cS\cU_C(2,2K_C)$. The map $\Phi_{K_C}$ is the composition of the theta map $\theta_{K_C}:\cS\cU(3,3K_C)\ra |3\Theta|$  with  the linear   projection ${\pi}_e \colon \vert 3 \Theta \vert \dashrightarrow \vert 3K_C \vert$ whose  center is $L_e \simeq {\PP}^3$. 
\end{proposition}

\begin{proof}
Recall from section \ref{secdet} that $\Phi_{K_C}$ factors through ${\pi}_e$. Since $dim \ |3K_C| \ =4$ this implies that the center of the projection  is a 3 dimensional linear subspace. 
By Remark \ref{g2} we know that the decomposable bundles of type $E\oplus K_C$, $E \in \cS\cU_C(2,2K_C)$, have $h^0>3$. This locus in $\cS\cU(3,3K_C)$ has dimension 3, and it is contained in the indeterminacy locus of ${\phi}_{K_C}$.
Now, the hyperelliptic involution $h$ on $C$ defines the involution 

\begin{eqnarray*}
\sigma: \cS\cU(3,\cO_C) & \ra &  \cS\cU(3,\cO_C);\\
F & \mapsto & h^*F^*;
\end{eqnarray*}

which is associated  to the 2:1 covering given by $\theta$ \cite{Quang}. In \cite{Quang} (Sect. 3 and 4),  it is shown that the  locus  given by vector bundles $F = E \oplus \cO_C$,  $E\in \cS\cU_C(2,\cO)$,  is contained in the fixed locus of $\sigma$  and it is embedded in $\cS\cU(3,\cO_C)$. The twist by $K_C$ of these
 vector bundles give those of type $E\oplus K_C$ with $E \in \cS\cU(2,2K_C)$. The image via $\theta$ of this locus is  $L_e \simeq {\PP^3}$ and  it is actually the center of the projection.
\end{proof}


Now let us recall some results from the literature about $\cS\cU(3,\cO_C)$.
Of course the twist by $K_C$ is an isomorphism and it is easy to understand which is the analogue result for $\cS\cU_C(3,3K_C)$. Since $\PP^8$ is smooth, the image of the singular locus of $\cS\cU(3,\cO_C)$ and the singular locus of the branch locus coincide, i.e. $Sing(\cC_6)=\theta(Sing(\cS\cU(3)))$. On $\pic^1(C)$ we have the involution $\lambda: L \mapsto K_C \otimes L^{-1}$ that leaves $\Theta$ invariant.
%
%
Hence $\lambda$ induces an action on all the powers of $\Theta$ and in particular, on $|3\Theta|$. The linear system $|3\Theta|$ decomposes in two eigenspaces, respectively 4 and 3 dimensional. We call $\PP^4_e$ the 4-dimensional eigenspace. It turns out that it  cuts out on $\cC_6$ a reducible variety given by a double $\PP^3$ (which is indeed contained in $Sing(\cC_6)$) and a quartic hypersurface $I \subset \PP^4$.  After the twist by $K_C$,  the first component is precisely  $L_e$, whereas the quartic 3-fold is an Igusa quartic. 

\begin{lemma}
The intersection of the closure of the general fiber of $\pi_e$ with $Sing(\cC_6)$ is a $15_3$ configuration of lines and points.
\end{lemma}

\begin{proof}
Recall that $Sing(\cC_6)$ is the locus of theta divisors corresponding to decomposable bundles, 
we will prove the claim by constructing explicitly these bundles.
Let us denote $\Delta_{3K_C}$ the closed subset of $|3K_C|$ given by the intersection with the big diagonal of $C^{(6)}$. 
Let us take $G=q_1+ \dots +q_6\in |3K_C| - \Delta_{3K_C}$, and let us consider the fiber of $\Phi_{K_C}$ over $G$. In order to guarantee the semi-stability of the vector bundles, the only totally decomposable bundles in the fiber of $G$ are all the 15 obtained by permuting the $q_i$'s in $\cO_C(q_1 + q_2)\oplus \cO_C(q_3+q_4)\oplus \cO_C(q_5+q_6)$. Let us now consider the bundles that decompose as the direct sum of a line bundle $L$ and a rank two indecomposable bundle. By the previous argument of semi-stability then $L$ must be of the type $\cO_C(p_i+p_j)$ for some $i,j \in \{1,\dots,6\}$ and $E$ should have fundamental divisor $\DD_E=\sum_{k\neq i,j} p_k$. Call $F$ the line bundle $\cO_C(\sum_{k\neq i,j} p_k) \equiv 3K_C -p_i-p_j$. It is easy to see that $\cS\cU_C(2, F)\cong \cS\cU_C(2,\cO_C)\cong \PP^3$, the isomorphism being given by the tensor product by a square root $F'$ of $F$. Now we recall from \cite{bol:kumwed} the following description of the fundamental map $\Phi_{F'}:\cS\cU_C(2, F) \dashrightarrow |F|$. The linear system $|F|$ is a $\PP^2$ and the fibers of $\Phi_{F'}$ are just lines  passing by $D\in |F|$ and the origin $[\cO_C\oplus\cO_C]$. Now the composition of the following embedding

\begin{eqnarray*}
\zeta: \cS\cU_C(2, F) & \hookrightarrow & \cS\cU_C(3,3K),\\
E & \mapsto & \cO_C(p_i+p_j) \oplus E,
\end{eqnarray*}
 
with the theta map is linear. In fact $\zeta(\cS\cU_C(2, F))$ is contained in the branch locus and the associated $3\Theta$ divisors form a three dimensional linear subsystem isomorphic to $|2\Theta|$. 
Then the image of $\zeta$ intersects the closure $\overline{\Phi^{-1}_{K_C}(G)}$ of the fiber over $G$ exactly along the fiber of $\Phi_{F'}$ over the divisor $\sum_{k\neq i,j} p_k \in |F|$. By \cite{bol:kumwed} we know that this is a line and it is not difficult to see that it contains 3 of the 15 totally decomposable bundles. On the other hand each totally decomposable bundle with fundamental divisor $G$ is contained in three lines of this kind.

\end{proof}

\begin{remark}\label{15dege}
When the divisor $G$ is taken in $\Delta_{3K_C}$ then the configuration $15_3$ degenerates because some of the points and of the lines coincide.
\end{remark}

\begin{theorem}\label{fibrusa}
The closure of the general fiber of ${\Phi}_{K_C}$ is the GIT quotient \\ $(\PP^2)^6//PGL(3)$.
\end{theorem}

\begin{proof}
We recall that $L_e$ is contained in $Sing(\cC_6)$ and in particular \newline scheme-theoretically it is contained twice in $\cC_6$. Since ${\Phi}_{K_C} $ factors through the projection with center $L_e$, then $\overline{\Phi_{K_C}^{-1}(G)}$, for $G \in |3K_C| - \Delta_{3K_C}$, is a degree two cover of $\PP^4_G:= \overline{\pi_e^{-1}(G)}$ ramified along the intersection of $\cC_6$ with $\PP^4_G$ which is residual to $2L_e$. This intersection is then a quartic hypersurface in $\PP^4_G$. Notably, since $L_e \subset \PP_e^4$, there exist a point $T \in |3K_C|$ s.t. $\overline{\pi_e^{-1}(T)} \cap \cC_6 $ is an Igusa quartic (see Prop. 5.2 of \cite{oxpa:prv} or \cite{minhigusa} Sect. 4).


Let us now blow up $|3\Theta|$ along $L_e$ and call $\widetilde{\PP}^8$ the obtained variety, which contains canonically the blown up Coble sextic, that we denote $\widetilde{C}_6$. Then the rational map $\pi_e$ resolves in a proper, flat map $\widetilde{\pi}_e$ as in the following diagram.

\qquad \qquad \qquad \qquad \qquad $\xymatrix{ \widetilde{\PP}^8 \ar[d] \ar[dr]^{\widetilde{\pi}_e} &  \\
|3\Theta| \ar@{-->}[r]^{{\pi}_e}& |3K_C| }$

This boils down to saying that the restriction of $\widetilde{\pi}_e$ to $\widetilde{C}_6$ is a flat family of quartic 3-folds over $|3K_C|$ and for any $B\in |3K_C|$ we have an isomorphism $\overline{\pi^{-1}_e(B)}\cap \cC_6 \cong \widetilde{\pi}^{-1}_e(B)\cap \tilde{C}_6$. Hence also one fiber of $\widetilde{\pi}_{e|\widetilde{C}_6}$ is an Igusa quartic. Since the ideal of the singular locus of $\cI_4$ is generated by the four polar cubics (\cite{hu:gsq}, Lemma 3.3.13) then the Igusa quartic has no infinitesimal deformations, i.e. it is rigid. This implies that the generic member of the flat family of quartics over $|3K_C|$ is an Igusa quartic. This in turn implies that the closure of the generic fiber of $\Phi_{K_C}$ is $(\PP^2)^6//PGL(3)$.
\end{proof}




\begin{corollary}
The Coble sextic $\cC_6$ is birational to a fibration over $\PP^4$ whose fibers are Igusa quartics.
\end{corollary}

\begin{corollary}
Along the generic fiber of $\Phi_{K_C}$, the involution of the degree two covering $\cS\cU(3,\cO_C)\ra |3\Theta|$ coincides with the involution given by the association isomorphism on $(\PP^2)^6//PGL(3)$.
\end{corollary}


\medskip

The fact that the intersection of $Sing(\cC_6)$ with the fibers over the open set $|3K_C| - \Delta_{3K_C}$ is precisely a $15_3$ configuration makes us argue that $|3K_C| - \Delta_{3K_C}$ should be the open locus where by rigidity the family of quartic three-folds is isotrivially isomorphic to the Igusa quartic. As already seen in Remark \ref{15dege}, if $B$ is an effective divisor out of this locus then, $\overline{\pi_e^{-1}(B)}\cap Sing(\cC_6)$ is a \textit{degenerate} $15_3$ configuration, in the sense that some of the 15 points and lines come to coincide. We are not able to prove the following, but it is tempting to say that this is all the singular locus of the special quartic three-folds over $\Delta_{3K}$. These would give very interesting examples of \textit{degenerate Igusa quartics}.
It would be interesting to study projective properties of these fibers such as the relation with the Segre cubics or with the Mumford-Knudsen compactification $\overline{\cM}_{0,6}$ of the moduli space of 6 points on a line. For instance, do they come from linear systems on $\overline{\cM}_{0,6}$? If it is so, what linear systems on $\overline{\cM}_{0,6}$ do they come from?


\medskip

The rational dual map of the Coble sextic has been thoroughly studied and described in \cite{orte:cob} and \cite{Quang}. Let us denote by $X_0,\dots ,X_8$ the coordinates on $\PP^8\cong |3\Theta|$ and by $F(X_0:\dots:X_8)$ the degree six poynomial defining $\cC_6$. Then the dual map is defined as follows:

\begin{eqnarray*}
\cD_6:\cC_6 & \dashrightarrow & \cC_3;\\
x & \mapsto & [\frac{\partial F}{\partial X_0}: \dots : \frac{\partial F}{\partial X_8}].
\end{eqnarray*}

The polar linear system is given by quintics that vanish along $Sing(\cC_6)$. Now fix a general divisor $B\in |3K|$ and call  $\cI_B$ the Igusa quartic defined by $(\cC_6\cap \overline{\pi_e^{-1}(B)}) - 2L_e$. Let us consider the restriction of $\cD_6$ to $\cI_B\subset \overline{{\pi}_e^{-1}(B)} =: \PP^4_B$ and denote by $A$ the $15_3$ configuration of points and lines in  $\PP^4_B$. Let now $H$ be the class of $L_e$ in $\pic(\PP^4_B)$ and consider the $4$- dimensional linear system $|\cI_A(3)+2H|$ on $\cI_B$. We can show the following.

\begin{proposition}\label{igudual}
The restricted dual map $\cD_{6|\cI_B}$ is given by a linear system $|\cD_{\cI_B}|$ that contains $|\cI_A(3)+2H|$ as a linear subsystem.
\end{proposition}


Remark that this means that for the general fiber $\cI_B$, there exists a canonical way to project the image $\cD(\cI_B)\subset \cC_3$ to a $\PP^4_B$ where the image of $\cI_B$ is a Segre cubic. This is resumed in the following.

\medskip

\begin{corollary}
The Coble cubic is birational to a fibration in Segre cubics over  $\PP^4$.
\end{corollary}

\begin{remark}
The birationality in itself is trivial, since $\cC_3$ is birational to $\cC_6$ which is birational to a fibration in Igusa quartics (which are in turn all birational to the Segre cubic) over $\PP^4$. The projections on the linear systems $|\cI_A(3)+2H|$ give a constructive canonical way to realize it. 
\end{remark}



\subsection{The Coble quartic}

In his subsection we assume that $C$ is a curve of genus $3$ and we consider the  moduli space $\cS\cU(2,\cO_C)$.  We recall that  the Kummer variety $Kum(C):=J(C)/\pm Id$ of $C$ is contained naturally in the $2\Theta$-linear series, whereas the moduli space $\cS\cU(2,\cO_C)$ is embedded by $\theta$ in $\PP^7\cong |2\Theta|$ as the unique quartic hypersurface $\cC_4$ singular along $Kum(C)$. This hypersurface is called the \textit{Coble quartic}. It is also known \cite{paulydual} that the Coble quartic is projectively self dual.
\hfill\par
Now we need to introduce a second important modular variety, i.e. the \textit{Segre cubic}. This is a nodal (and hence rational) cubic three-fold $S_3$ in $\PP^4$ whose singular locus is given by ten double points. There is a natural action of $\Sigma_6$ on this projective space and $S_3$ is invariant with respect to this action. $S_3$ is in fact the GIT quotient $(\PP^1)^6//PGL(2)$ \cite{do:pstf}. Moreover, $S_3$ realizes the so-called \textit{Varchenko bound}, that is, it has the maximum number of double points (ten) that a cubic threefold with isolated singularities may have and  this property identifies the Segre cubic in a unique way, up to projective equivalence. As already stated it is the projective dual variety of the Igusa quartic.

Our construction allows us to give a simple proof of the following result from \cite{albol}.

\begin{proposition}\label{quartic}
The moduli space $\cS\cU(2,\cO_C)$ is birational to a fibration over  $\PP^3$ whose fibers are Segre cubics.
\end{proposition}

\begin{proof}
After the customary twist by a degree 3 divisor $D$, the fundamental map is ${\Phi}_D \colon  \cS\cU(2,\cO_C(2D)) \to \vert \cO_C(2D) \vert\cong \PP^3.$ Since $\theta_D$ is an embedding, we identify  $\cS\cU(2,\cO_C(2D))$ with its image $\cC_4\subset \PP^7$ and  ${\Phi}_D$ with the linear projection  onto  $|2D|$. The center of the projection is the linear span of the locus of vector bundles $E$ s.t. $h^0(C,E\otimes\cO_C(D))>2$. If $E$ is stable then $E\cong E^*$ and by an easy Riemann-Roch computation we find that $h^0(C,E\otimes\cO_C(D))>2$ if and only if $h^0(C,E\otimes \cO_C(K-D)) > 0$. As it is shown in \cite{albol} Prop 3.1 this is equivalent to the fact the $E$ lies in the $\PP^3\cong|3K-2D|^*\subset \cC_4$ that parametrizes vector bundles $E$ that can be written as an extension of the following type

$$0\ra \cO_C(D-K) \ra E \ra \cO_C(K-D)\ra 0.$$

Let us denote $\PP^3_c$ this projective space. $\cC_4$ contains $\PP^3_c$ with multiplicity one. This implies that the closure of any fiber of the projection ${\Phi}_D :\cC_4 \dashrightarrow \vert 2 D \vert$ is a cubic 3-fold contained in the $\PP^4$ spanned by $\PP^3_c$ and a point of $\vert 2 D \vert $. Let us denote as usual $\Delta_D$ the intersection of the large diagonal with the linear system $|2D|\subset C^{(6)}$. Then suppose we fix a $B\in |2D| - \Delta_D$. Let us consider the intersection of the fiber of $\Phi_D$ over $B$ with the strictly semi-stable locus. By semi-stability it is easy to see that these points correspond to the partitions of the 6 points of $B$ in complementary subsets of 3 elements each. We have ten of them. As stated her above, a cubic 3-fold can not have more than ten ordinary double points and the Segre cubic is uniquely defined by this singular locus up to projective equivalence.
\end{proof}


Also in this case, if $B\in \Delta_D$ then the intersection $\overline{\Phi_D^{-1}(B)}\cap Kum(C)$ is set-theoretically a finite set of points of cardinality strictly smaller then 10: the singular locus seems to degenerate. It is tempting, like in the case of Igusa quartics, to say that some of these points have multiplicity bigger than one and we obtain \textit{degenerate Segre cubics} over $\Delta_D$.

\medskip

As we have already remarked, also in the case of $\cC_4$ the polar map is well known and described. Let $Y_i$ be the coordinates on $\PP^7\cong |2\Theta|$ and $G(Y_0:\dots:Y_7)$ the quartic equation defining $\cC_4$, then the (self) polar rational map of $\cC_4$ is defined in the following way.

\begin{eqnarray*}
\cD_4:\cC_4 & \dashrightarrow & \cC_4;\\
x & \mapsto & [\frac{\partial G}{\partial Y_0}: \dots : \frac{\partial G}{\partial Y_7}].
\end{eqnarray*}

Let $B\in |2D| - \Delta_D$ and let $\PP^4_B$ be the linear span of the point corresponding to $B$ and of $\PP^3_c$
It turns out that the restriction of $\cD_4$ to $\PP^4_B$ beahaves in a way very similar to the case of the sextic (see Prop. \ref{igudual}).. Let $S_{3B}\subset \PP^4_B$ be the Segre cubic such that $\cC_4 \cap \PP^4_B = S_{3B} \cup \PP^3_c$. We denote $J$ the set of 10 nodes of $S_{3B}$. Then the linear series $|\cI_J(2)|$ on $S_{3B}$ is the polar system of the Segre cubic.

\begin{proposition}\label{segdual}
The restricted dual map $\cD_{4|S_{3B}}$ is given by a linear system $|\cD_{S_3}|$ that contains $|\cI_J(2)+H|$ as a linear subsystem.
\end{proposition}

As in the case of $\cC_6$ this implies that we have a canonical way to construct the birational map of the following corollary via the polar map $\cD_4$.

\begin{corollary}
The Coble quartic is birational to a fibration in Igusa quartics over $\PP^3$
\end{corollary}



\bibliographystyle{amsalpha}
\bibliography{bibkaji}

 \end{document}